\documentclass[12pt]{article}

\usepackage{amsmath,amssymb,url}
\usepackage{tikz,color,pgf}
\usepackage{longtable} 
\usepackage{enumitem}

\newtheorem{theorem}{Theorem}[section]

\newtheorem{lemma}[theorem]{Lemma}

\newcommand{\proof}{\noindent{\bf Proof.\ }}
\newcommand{\qed}{\hfill $\square$ \bigskip}

\begin{document}

\title{Correcting the algorithm for the secure domination number of cographs by Jha, Pradhan, and Banerjee}

\author{Anja Ki\v sek$^{a}$\thanks{Email: \texttt{anja.kisek@gmail.com}} 
\and Sandi Klav\v zar $^{a,b,c}$\thanks{Email: \texttt{sandi.klavzar@fmf.uni-lj.si}}
}
\maketitle

\begin{center}
$^a$ Faculty of Mathematics and Physics, University of Ljubljana, Slovenia\\
\medskip

$^b$ Institute of Mathematics, Physics and Mechanics, Ljubljana, Slovenia\\
\medskip

$^c$ Faculty of Natural Sciences and Mathematics, University of Maribor, Slovenia\\
\medskip
\end{center}

\begin{abstract}
Jha, Pradhan, and Banerjee devised a linear algorithm to compute the secure domination number of a cograph. Here it is shown that their Lemma~2, which is crucial for the computational complexity of the algorithm, is incomplete. An accordingly modified lemma is proved and it is demonstrated that the complexity of the modified algorithm remains linear. 
\end{abstract}


\noindent
{\bf AMS Subj.\ Class.\ (2020)}: 05C57, 05C69, 05C76, 05C85

\section{Introduction and preliminaries}
\label{sec:intro}

Cographs have been independently introduced in many different contexts and under many different names (including D*-graphs, hereditary Dacey graphs, $2$-parity graphs), clearly indicating the intrinsic role of this class of graphs in graph theory and elsewhere. Cographs can be characterized in several different ways, see~\cite{corneil-1981}. Moreover, their research is still very active, the list~\cite{allem-2020, bresar-2015, epple-2021, geis-2020, ghorbani-2019, tsur-2020} presents just a sample of recent studies and applications of cographs. 

Numerous problems, including the independence number, the clique number, and the existence of a Hamiltonian cycle, that are difficult in general, are polynomial (in most cases linear) on cographs. In this direction of research, two linear algorithms have been independently designed that compute the secure domination number of a cograph~\cite{araki-2019, jha-2019}.   

The secure domination problem was introduced in~\cite{cockayne-2003}, see also~\cite{cockayne-2005}, motivated with a problem of assigning guards at various locations corresponding to a dominating set. Computing the secure domination number is NP-complete in general, and remains NP complete  when restricted to bipartite and split graphs~\cite{boumediene-2015}, star convex bipartite graphs and doubly chordal graphs~\cite{wang-2018}, and chordal bipartite graphs~\cite{pradhan-2018}. On the positive side, the problem is linear on trees~\cite{burger-2014} (cf.\ also~\cite{li-2017}) and more generally on block graphs~\cite{pradhan-2018}, as well as on proper interval graphs~\cite{araki-2018}. We also point to~\cite{burdett-2020+}, where the secure domination problem is approached via a binary programming formulation of the problem. 

We proceed as follows. In the rest of this section definitions needed are listed. In Section~\ref{sec:lemma-original} we present the necessary set up for the algorithm from~\cite{jha-2019}, and construct an infinite family of cographs which demonstrate that~\cite[Lemma~2]{jha-2019} does not hold in general. The latter lemma is a key ingredient for the linearity of the designed algorithm. In Section~\ref{sec:lemma-corrected} we then propose a modified lemma, prove it, and argue that the algorithm remains linear also when the new lemma is applied. 

We now proceed with definitions, for other basic graph theory concepts not defined here we follow~\cite{west-2001}. Let $G = (V(G), E(G))$ be a graph. A \emph{clique} of $G$ in a complete sugraph of $G$. By abuse of language we will also say that the set of vertices which induces a complete graph is a clique. The \emph{open neighbourhood} of $v$ in $G$ is the set $N_G(v) = \{u:\ uv \in E(G)\}$, the \emph{closed neighbourhood} of $v$ is $N_G[v] = N_G(v) \cup \{v\}$. The \emph{complement} $\overline{G}$ of $G$ has $V(\overline{G}) = V(G)$ and $E(\overline{G}) = \{(u,v):\ (u,v)\notin E(G) \}$. If $G_1$ and $G_2$ are disjoint graphs, then their \emph{union} is the graph with the vertex set $V(G_1) \cup V(G_2)$ and the edge set $E(G_1) \cup E(G_2)$. The \emph{join} $G_1 + G_2$ of $G_1$ and $G_2$ has $V(G_1 + G_2) = V(G_1) \cup V(G_2)$ and $E(G_1 + G_2) = E(G_1) \cup E(G_2) \cup \{(u,v):\ u \in V(G_1), v \in V(G_2) \}$. Note that $G_1 + G_2 = \overline{\overline{G_1} \cup \overline{G_2}}$. 

$D\subseteq V(G)$ is a {\em dominating set} of a graph $G$ if each vertex $x \in V(G) \setminus D$ is adjacent to a vertex from $D$. The \emph{domination number} $\gamma(G)$ of $G$ is the cardinality of a smallest dominating set of $G$. A dominating set $S\subseteq V(G)$ is a {\em secure dominating set} of $G$  if for each vertex $x\in V(G)\setminus S$ there exists its neighbour $y\in S$ such that $(S \cup\{x\})\setminus \{y\}$ is a dominating set. The \emph{secure domination number} $\gamma_s(G)$ of $G$  is the minimum cardinality of a secure dominating set of $G$.

\emph{Cographs} are defined resursively as follows: (i) $K_1$ is cograph, (ii) if $G$ is a cograph, then $\overline{G}$ is a cograph, and (iii) if $G_1, \dots, G_k$, $k \geq 2$, are cographs, then $G_1 \cup \cdots \cup G_k$ is a cograph. A cograph $G$ can be represented as a rooted tree $T_G$ called a \emph{cotree} of $G$~\cite{corneil-1981}.  In this representation, leaves of $T_G$ are vertices of the original cograph $G$. Inner vertices of $T_G$ are labeled by $\cup$ or by $+$, depending, respectively, whether the cotree rooted in the inner vertex corresponds to a cograph, obtained by the union or the join of the cographs associated with its children. (See Fig.~\ref{fig:counterExampleTree} for the cotree of the cograph from Fig.~\ref{fig:counterExampleGraph}.) Finally, if $v\in T_G$, then we denote by $T_G(v)$ the subgraph of $G$ induced by the leaves of the subtree of $T_G$ rooted at the vertex $v$.

\section{Original lemma and counterexamples}
\label{sec:lemma-original}

Let $T_G$ be the cotree of a cograph $G$ and let $c$ be an inner vertex of $T_G$. If $c$ has label $\cup$, then each children of $c$ is either a leaf or a vertex with label $+$, and the cograph $T_G(c)$ is disconnected. If $c$ has label $+$, then each children of $c$ is either a leaf or a vertex with label $\cup$.

Following~\cite{jha-2019}, we assign a label $\mathcal{R}$ to a vertex $u$ of $T_G$ if (i) $u$ has label $\cup$, (ii) $u$ has exactly two children in $T_G$, say $x$ and $y$, (iii) $\gamma(T_G(x)) = 1$, and (iv) $\gamma_s(T_G(y)) = 1$. The following lemma gives a simple characterization of the vertices of $T_G$ to which label $\mathcal{R}$ is assigned. 

\begin{lemma}\label{Lemma1} {\rm \cite[Lemma~1]{jha-2019}}
Let $T_G$ be the cotree of a cograph $G$. Then $u\in V(T_G)$ has label $\mathcal{R}$ if and only if $T_G(u)$ is disconnected and there exists a vertex $w\in V(T_G(u))$ such that $V(T_G(u)) \setminus N_{T_G(u)}[w]$ is a clique.
\end{lemma}

We next recall when a cograph has property $\mathcal{P}$ which is a key for a fast computation of the secure domination number of a cograph.  Let $G$ be a cograph that is the join of cographs $G_1, \dots, G_\ell$, $\ell \geq 2$. Then we say that $G$ satisfies \emph{property $\mathcal{P}$} if there exist two distinct vertices $x,y \in V(G)$ such that $\{x, y\}$ is a dominating set of $G$ and each of $V(G) \setminus N_G[x]$ and $V(G) \setminus N_G[y]$ is either empty or a clique. The following characterization of the cographs with property $\mathcal{P}$ was claimed. 

\begin{lemma}\label{Lemma2Original} {\rm \cite[Lemma~2]{jha-2019}} 
Let $T_G$ be the cotree of a cograph $G$ and let $c$ be a vertex of $T_G$ with label $+$. Then the cograph $T_G(c)$ satisfies property $\mathcal{P}$ if and only if there exist at least two children of $c$ in $T_G$ such that each of them is either a leaf or a vertex with label $\mathcal{R}$.
\end{lemma}

We now provide an infinite family of cographs which demonstrate that Lemma~\ref{Lemma2Original} is not true as stated. Let $k \in \mathbb{N}$, and let $K_k$ be a complete graph with $V(K_k) = \{a_1, \dots, a_k\}$. Define $G_k$ to be the cograph with the vertex set $V(G_k) = V(K_k) \cup \{b, c, d, e\}$ and the edge set 
$$E(G_k) = E(K_k) \cup \{bc, bd, be\} \cup \{ca_i, da_i, ea_i:\  i\in [k]\}\,.$$

\begin{figure}[t!]
\centering
\begin{tikzpicture}[main_node/.style={circle,draw,minimum size=2em,inner sep=1, scale=0.9}]
\usetikzlibrary{shapes}
\pgfdeclarepatternformonly{mynewdots}{\pgfqpoint{-1pt}{-1pt}}{\pgfqpoint{0.5pt}{0.5pt}}{\pgfqpoint{9pt}{9pt}}
{%
  \pgfpathcircle{\pgfqpoint{0pt}{0pt}}{.5pt}%
  \pgfusepath{fill}%
}%
\usetikzlibrary{patterns}
    \node[main_node] (1) at (-2,0) {$b$};
    \node[main_node] (2) at (0,1) {$c$};
    \node[main_node] (3) at (0,0) {$d$};
    \node[main_node] (4) at (0,-1) {$e$};
    \node[main_node, font=\bfseries] (5) at (2.5,0) {$a_i$};
    \draw[densely dotted, pattern=mynewdots] (2.5,0) ellipse (1.2 and 1.5);
    
    \draw (2) -- (1) -- (3);
    \draw (4) -- (1);
    \draw[thick] (2) -- (5) -- (3);
    \draw[thick] (5) -- (4);

\end{tikzpicture}
\caption{Graphs $G_k$. The dotted area represents the complete graph $K_k$, thicker edges indicate the join between $K_k$ and the independent set $\{c, d, e\}$.}
\label{fig:counterExampleGraph}
\end{figure}
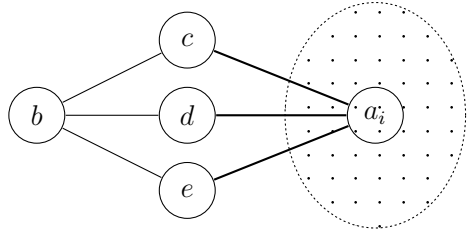

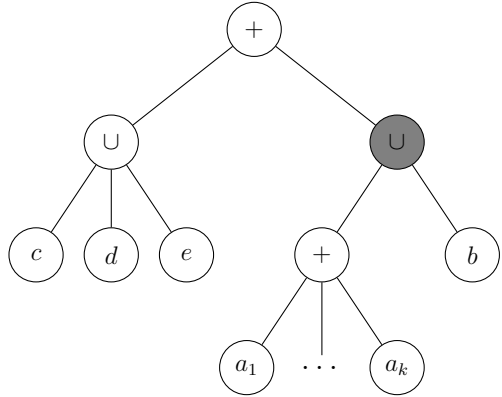
\begin{figure}[t!]
\centering
\begin{tikzpicture}[main_node/.style={circle,draw,inner sep=3pt,minimum size = 0.9cm, scale=0.8}]

	\node[main_node] (1)   at (1.1,0.5) {$+$};
    \node[main_node] (2)   at (-0.8, -1) {$\cup$};
	\node[main_node, fill=gray] (3)   at (3,-1){$\cup$};
    \node[main_node] (4)   at (-1.8, -2.5){$c$};
    \node[main_node] (5)   at (-0.8,-2.5) {$d$};
    \node[main_node] (6)   at (0.2,-2.5) {$e$};
    \node[main_node] (7)   at (2,-2.5) {$+$};
    \node[main_node] (8)   at (4,-2.5) {$b$};
    \node[main_node] (9)   at (1,-4) {$a_1$};
    \node[] (10)   at (2,-4) {$\dots$};
     \node[main_node] (11)   at (3,-4) {$a_k$};
    
    \draw (4) -- (2) -- (5);
    \draw (6) -- (2) -- (1) -- (3) -- (7) -- (9);
    \draw (8) -- (3);
    \draw (10) -- (7) -- (11);

\end{tikzpicture}
\caption{The cotree $T_{G_k}$. The gray vertex is labeled by $\mathcal{R}$.}
\label{fig:counterExampleTree}
\end{figure}

The graphs $G_k$, $k \in \mathbb{N}$, are schematically shown  in~Fig.~\ref{fig:counterExampleGraph}, while their cotrees $T_{G_k}$ are drawn in Fig.~\ref{fig:counterExampleTree}. 

Let $x$ be the root of $T_{G_k}$, its label is thus $+$. We first claim that $T_{G_k}(x) = G_k$ satisfies property $\mathcal{P}$. To see it observe that $\{a_1, b\}$ is a dominating set, and that $V(G) \setminus N_G[a_1]=\{b\}$ and $V(G) \setminus N_G[b] = V(K_k)$ are both cliques. We next claim that $G_k$ does not fulfill the condition of Lemma~\ref{Lemma2Original}. The vertex $x$ has in $T_{G_k}$ two children labelled with $\cup$. It is easy to check that the right child (the gray vertex in Fig.~\ref{fig:counterExampleTree}) is labelled with $\mathcal{R}$, but the left child (the one that is a parent of the leaves $c$, $d$ and $e$) is obviously not labelled with $\mathcal{R}$ due to the number of its children. Hence the left children of $x$ is neither a leaf nor a vertex with label $\mathcal{R}$, and therefore  the condition of Lemma~\ref{Lemma2Original} is not fulfilled. 

\section{Corrected lemma}
\label{sec:lemma-corrected}

In order to correct Lemma~\ref{Lemma2Original}, we give the following characterization of cographs which fulfill property $\mathcal{P}$. 

\begin{lemma}\label{Lemma2}
Let $T_G$ be the cotree of a cograph $G$ and let $c$ be a vertex of $T_G$ with label $+$. Then the cograph $T_G(c)$ satisfies property $\mathcal{P}$ if and only if there exist at least two children of $c$ in $T_G$ such that
\begin{enumerate}[label=(\roman*)]
\item each of them is either a leaf or a vertex with label $\mathcal{R}$, or
\item one of them is a vertex $u$ with label $\mathcal{R}$ and $T_G(u)$ is a union of two complete graphs.
\end{enumerate}
\end{lemma}

\proof
($\Rightarrow$) Assume that $T_G(c)$ has property $\mathcal{P}$. Then there exist two distinct vertices $x, y \in V(T_G(c))$ such that $\{x, y\}$ is a dominating set and each of $V(T_G(c)) \setminus N_{T_G(c)}[x]$ and $V(T_G(c)) \setminus N_{T_G(c)}[y]$ is either empty or a clique. By definition of $\mathcal{P}$, the cograph $T_G(c)$ is the join of at least two cographs, hence at least two children of $c$ are guaranteed.

If $V(T_G(c)) \setminus N_{T_G(c)}[x] =  \emptyset$, then $x$ is adjacent to each vertex of $V(T_G(c)) \setminus \{x\}$. As $c$ has label $+$, the vertex $x$ must be a leaf adjacent to $c$, for otherwise there would exist some vertex of $T_G(c)$ not adjacent to $x$.

If $V(T_G(c)) \setminus N_{T_G(c)}[x]$ is a clique, then with a parallel reasoning we conclude that $x$ is not a child of $c$. Let $w_x\in V(T_G(c))$ be a child of $c$, such that $T_G(w_x)$ contains $x$. By the observations from the first paragraph of Section~\ref{sec:lemma-original} we know that $w_x$ is labeled by $\cup$ and that $T_G(w_x)$ is disconnected. Since $x \in V(T_G(w_x))$ and $c$ has label $+$, the vertex $x$ is adjacent to each vertex of $V(T_G(c)) \setminus V(T_G(w_x))$. Moreover, because $V(T_G(c)) \setminus N_{T_G(c)}[x]$ is a clique, also $V(T_G(w_x)) \setminus N_{T_G(w_x)}[x]$ is a clique. Now we have a disconnected graph $T_G(w_x)$ containing $x$, such that $V(T_G(w_x)) \setminus N_{T_G(w_x)}[x]$ is a clique. By Lemma~\ref{Lemma1}, $w_x$ has label $\mathcal{R}$.

Similarly, if $V(T_G(c)) \setminus  N_{T_G(c)}[y] = \emptyset$, then $y$ is a child of $c$ such that $y$ is a leaf. Otherwise there exists a child of $c$, say $w_y$, such that $y \in V(T_G(w_y))$ and $w_y$ has label $\mathcal{R}$. In case that both of $V(T_G(c)) \setminus N_{T_G(c)}[x]$ and  $V(T_G(c)) \setminus N_{T_G(c)}[y]$ are cliques, we have two vertices $w_x$ and $w_y$ such that they are children of $c$ and have label $\mathcal{R}$. If $w_x \neq w_y$, vertex $c$ has at least two children, both labeled with $\mathcal{R}$.

In the cases discussed so far, the vertex $c$ has at least two children and each of them is either a leaf or a vertex with label $\mathcal{R}$, hence (i) holds. It can happen, however, that $w_x=w_y=w$ holds. In this case $x, y \in T_G(w)$. Since $w$ has label $\mathcal{R}$, the graph $T_G(w)$ is disconnected and has two components, $H_1$ and $H_2$. The assumption that $\{x, y\}$ is a dominating set implies that, without loss of generality, $x \in H_1$ and $y \in H_2$. Since $V(T_G(c)) \setminus N_{T_G(c)}[x] = V(H_2)$, the subgraph $H_2$ must be a clique. Similarly, since $V(T_G(c)) \setminus N_{T_G(c)}[y] = V(H_1)$, also $H_1$ is a clique. In this final case we thus see that (ii) holds. 

\medskip
($\Leftarrow$) 
Suppose first that there exist at least two children of $c$ in $T_G$ such that each of them is either a leaf or a vertex with label $\mathcal{R}$: Let $u,v$ be such children of $c$. First assume $u$ is a leaf. Since $u \in T_G(c)$ and $u$ is adjacent to each vertex of $T_G(c) \setminus \{u\}$, we infer that  $V(T_G(c)) \setminus N_{T_G(c)}[u] = \emptyset$. If $v$ is also a leaf, then also $V(T_G(c)) \setminus N_{T_G(c)}[v] = \emptyset$. Hence $T_G(c)$ satisfies property $\mathcal{P}$.

In case $v$ is not a leaf, it has label $\mathcal{R}$. By Lemma~\ref{Lemma1}, $T_G(v)$ is disconnected and there exists a vertex $w \in T_G(v)$ such that $V(T_G(v)) \setminus N_{T_G(v)}(w)$ is clique. Since $c$ has label $+$, the vertex $w$ is adjacent to each vertex in $T_G(c) \setminus T_G(v)$. Since  $T_G(v) \setminus N_{T_G(v)}(w)$ is a clique, $T_G(c) \setminus N_{T_G(v)}(w)$ is also a clique. Since $u$ is a leaf, the set $\{u, v\}$ is dominating, hence $T_G(c)$ satisfies property $\mathcal{P}$.

Now assume that both $u$ and $v$ have label $\mathcal{R}$. As discussed above, there exist $w_1 \in V(T_G(u))$ and $w_2 \in V(T_G(v))$ such that $V(T_G(c)) \setminus N_{T_G(c)}[w_1]$ and $V(T_G(c)) \setminus N_{T_G(c)}[w_2]$ are cliques. To satisfy property $\mathcal{P}$ we need to see that $\{w_1, w_2\}$ is dominating set. This is indeed the case since $c$ has label $+$, and therefore $w_1$ is adjacent to each vertex of $V(T_G(c)) \setminus V(T_G(u))$, and $w_2$ is adjacent to each vertex of  $V(T_G(c)) \setminus V(T_G(v))$.

\medskip
Suppose second that there exist at least two children of $c$ in $T_G$, and one of them is a vertex $u$ with label $\mathcal{R}$ and $T_G(u)$ is a union of two complete graphs. Let $x$ and $y$ be arbitrary vertices from each of the components of $T_G(u)$. Since $c$ is labeled by $+$, $x$ and $y$ are adjacent to each vertex in $G \setminus T_G(u)$. Since both components of $T_G(u)$ are cliques, both $V(G) \setminus N_{T_G(c)}[x] = N_{T_G(u)}[y]$ and $V(G) \setminus N_{T_G(c)}[y] = N_{T_G(u)}[x]$ are cliques. Also notice that $N_{T_G(c)}[x] \cup N_{T_G(c)}[y] = V(T_G(c))$, so $\{x, y\}$ is a dominating set. Hence $T_G(c)$ satisfies property $\mathcal{P}$.
\qed

As we have already mentioned, an efficient testing whether a given cograph has property $\mathcal{P}$ is crucial for the linearity of the algorithm from~\cite{jha-2019} which computes the secure domination number of cographs. More precisely, on each iteration of the algorithm, in case the current vertex of the cotree $T_G$ is labeled by $+$,  a checking of property $\mathcal{P}$ is required. Authors in~\cite{jha-2019} provided an implementation with an array, tracking labels $\mathcal{R}$ for each vertex of the cotree, and showed that it can be updated in $O(|N_{T_G}(c_i)|)$ time. The new Lemma~\ref{Lemma2} requires not only the data about the label $\mathcal{R}$ of vertices, but also completeness of subgraphs, hence a verification of property $\mathcal{P}$ needs to include values of $\gamma_s$ for descendants of the current vertex. Since the values $\gamma_s$ on each iteration of the algorithm are noted, algorithm remains linear also after changing verification of property $\mathcal{P}$ by means of Lemma~\ref{Lemma2}.

\section*{Acknowledgements}

We acknowledge the financial support from the Slovenian Research Agency (research core funding No.\ P1-0297 and projects J1-9109, J1-1693, N1-0095, N1-0108).


\begin{thebibliography}{99}

\bibitem{allem-2020}
  L.E.~Allem, F.~Tura, 
  Integral cographs,
  Discrete Appl.\ Math.\ 283 (2020) 153--167. 

\bibitem{araki-2018}
  T.~Araki, H.~Miyazaki, 
  Secure domination in proper interval graphs,
  Discrete Appl.\ Math.\ 247 (2018) 70--76. 

\bibitem{araki-2019}
  T.~Araki, R.~Yamanaka, 
  Secure domination in cographs, 
  Discrete Appl.\ Math.\ 262 (2019) 179--184.

\bibitem{boumediene-2015}
  H.Boumediene Merouane, M.~Chellali, 
  On secure domination in graphs,
  Inform.\ Process.\ Lett.\ 115 (2015) 786--790.

\bibitem{bresar-2015}
  B.~Bre\v sar, M.~Changat, T.~Gologranc, B.~Sukumaran, 
  Cographs which are cover-incomparability graphs of posets,
  Order 32 (2015) 179--187.
  
\bibitem{burdett-2020+}
  R.~Burdett, M.~Haythorpe, 
  An improved binary programming formulation for the secure domination problem,
  Ann.\ Oper.\ Res.\  295 (2020) 561--573.
  
\bibitem{burger-2014}
  A.P.~Burger, A.P.~de Villiers, J.H.~van Vuuren, 
  A linear algorithm for secure domination in trees,
  Discrete Appl.\ Math.\ 171 (2014) 15--27.

\bibitem{cockayne-2003}
  E.J.~Cockayne, O.~Favaron, C.M.~Mynhardt, 
  Secure domination, weak Roman domination and forbidden subgraphs,
  Bull.\ Inst.\ Combin.\ Appl.\ 39 (2003) 87--100. 

\bibitem{cockayne-2005}
  E.J.~Cockayne, P.J.P.~Grobler, W.R.~Grundlingh, J.~Munganga, J.H.~van Vuuren, 
  Protection of a graph,
  Util.\ Math.\ 67 (2005) 19--32.

 \bibitem{corneil-1981}
  D.G.~Corneil, H.~Lerchs, L.S.~Burlingham,
  Complement reducible graphs,
  Discrete Appl.\ Math.\ 3 (1981) 163--174.

\bibitem{epple-2021}
  D.A.~Epple, J.~Huang, 
  $(k,l)$-colourings and Ferrers diagram representations of cographs,
  European J.\ Combin. 91 (2021) Paper 103208.

\bibitem{geis-2020}
  M.~Gei\ss, P.F.~Stadler, M.~Hellmuth,
  Reciprocal best match graphs,
  J.\ Math.\ Bio.\ 80 (2020) 865--953. 
  
\bibitem{ghorbani-2019}
  E.~Ghorbani, 
  Cographs: eigenvalues and Dilworth number,
  Discrete Math.\ 342 (2019) 2797--2803.

\bibitem{jha-2019}
  A.~Jha, D.~Pradhan, S.~Banerjee, 
  The secure domination problem in cographs,
  Inform.\ Process.\ Lett.\ 145 (2019) 30--38.

\bibitem{li-2017}
  Z.~Li, Z.~Shao, J.~Xu, 
  On secure domination in trees,
  Quaest.\ Math.\ 40 (2017) 1--12. 

\bibitem{pradhan-2018}
  D.~Pradhan, A.~Jha, 
  On computing a minimum secure dominating set in block graphs,
  J.\ Comb.\ Optim.\ 35 (2018) 613--631.

\bibitem{tsur-2020}
  D.~Tsur, 
  Faster algorithms for cograph edge modification problems,
  Inform.\ Process.\ Lett.\ 158 (2020) Paper 105946. 

\bibitem{wang-2018}
  H.~Wang, Y.~Zhao, Y.~Deng, 
  The complexity of secure domination problem in graphs,
  Discuss.\ Math.\ Graph Theory 38 (2018) 385--396.

\bibitem{west-2001}
  D.B.~West,
  Introduction to Graph Theory, 2nd ed.,
  Prentice-Hall, NJ, 2001.
  
\end{thebibliography}
\end{document}